\documentclass[twoside,11pt]{amsart}

\usepackage{amssymb}
\usepackage{amsthm}

\newtheorem{theorem}{Theorem}[section]
\newtheorem{definition}[theorem]{Definition}
\newtheorem{lemma}[theorem]{Lemma}
\newtheorem{proposition}[theorem]{Proposition}
\newtheorem{question}[theorem]{Question}

\newcommand{\R}{\mathbb R}
\newcommand{\PP}{\mathcal P}
\newcommand{\HH}{\mathcal H}
\newcommand{\e}{\epsilon}
\newcommand{\QQ}{\mathcal Q}
\newcommand{\osc}{\mbox{osc} \,}
\newcommand{\leng}{\mbox{leng} \,}
\newcommand{\lenginf}{\mbox{leng}_\infty \,}
\newcommand{\lengoneinf}{\mbox{leng}_{(1,\infty)} \,}
\newcommand{\Dev}{\mbox{Dev} \,}

\title[THE GROUP OF HAMILTONIAN HOMEOMORPHISMS]{THE GROUP OF HAMILTONIAN \\ HOMEOMORPHISMS IN THE $L^\infty$-NORM}

\author{STEFAN M\"ULLER}

\email{mueller@math.wisc.edu}

\address{Department of Mathematics, University of Wisconsin-Madison, 480 Lincoln Drive, Madison WI 53706, USA}


\begin{document}

\thispagestyle{plain}

\begin{abstract}
\noindent The group $Hameo \, (M,\omega)$ of Hamiltonian homeomorphisms of a connected symplectic manifold $(M,\omega)$ 
was defined and studied in [7] and further in [6]. 
In these papers, the authors consistently used the $L^{(1,\infty)}$-Hofer norm 
(and not the $L^\infty$-Hofer norm) on the space of Hamiltonian paths (see below for the definitions). 
A justification for this choice was given in [7]. 
In this article we study the $L^\infty$-case. 
In view of the fact that the Hofer norm on the group $Ham \, (M,\omega)$ of Hamiltonian diffeomorphisms does not depend 
on the choice of the $L^{(1,\infty)}$-norm vs. the $L^\infty$-norm [9], 
Y.-G. Oh and D. McDuff (private communications) asked 
whether the two notions of Hamiltonian homeomorphisms arising from the different norms coincide. 
We will give an affirmative answer to this question in this paper. 
\end{abstract}

\maketitle


\tableofcontents


\section{Introduction}

Let $(M,\omega)$ denote a closed connected symplectic manifold. 
In section 2 we present a foundational study of the group of Hamiltonian homeomorphisms of $(M,\omega)$, 
which in the $L^{(1,\infty)}$-case appeared in [7] and [6]. 
See these papers for a more detailed study. 
Most of the theory goes through the same way, and in fact the proofs are often easier in the $L^\infty$-case, 
so we will frequently refer to these papers for details. 
We will compare the two norms and point out some of their advantages and disadvantages, 
again also see the papers cited above. 

Finally in section 3 we will give the proof of the main theorem, which states 

\begin{theorem} \label{firsttheorem}
$Hameo_\infty \, (M,\omega) = Hameo_{(1,\infty)} \, (M,\omega)$. \qed 
\end{theorem}

For the precise definitions of these objects, see section 2. 
Remark that the equality in the above theorem is between topological groups, 
where both sets are equipped with the $C^0$-topology (see section 2 for the definition) 
as subgroups of $Homeo \, (M)$, the group of homeomorphisms of $M$. 
In view of this theorem we may drop the subscripts. 

In the following, we will mainly use the notation from [7]. 
Denote by $C^\infty ([0,1] \times M)$ the vector space of 
smooth time-dependent Hamiltonian functions $H : [0,1] \times M \to \R$. 
Often we view $H$ as a family of smooth functions $H_t : M \to \R$. 
Each Hamiltonian $H \in C^\infty ([0,1] \times M)$ generates a family of diffeomorphisms $\phi_H^t$ of $M$, 
with $\phi_H^0 = id$, 
where $\phi_H^t (x)$ is the flow of the Hamiltonian vector field $X_H$ associated to the Hamiltonian $H$ defined by 
$$ \omega (X_H (t), \cdot) = d H_t, \hspace{5mm} \mbox{for all } t \in [0,1]. $$ 
That is, $\phi_H^t (x)$ solves Hamilton's equation $\dot x (t) = X_H (t, x(t) )$. 
We will always denote by $\phi_H$ the corresponding smooth Hamiltonian path $t \in [0,1] \mapsto \phi_H^t$, 
and by $\PP^{ham} (Symp \, (M,\omega), id)$ the set of all such Hamiltonian paths. 
Here $Symp \, (M,\omega)$ denotes the group of symplectic diffeomorphisms. 
The set of time-one maps of all such Hamiltonian paths is as usual denoted by $Ham \, (M,\omega)$, 
the group of Hamiltonian diffeomorphisms. 

Let $H \in C^\infty ([0,1] \times M)$. 
Denote by 
$$ \osc (H_t) = \max_{x \in M} H_t (x) - \min_{x \in M} H_t (x), \hspace{5mm} t \in [0,1], $$ 
the oscillation of $H_t$. 
Then define 
$$ \| H \|_{(1,\infty)} = \int_0^1 \osc (H_t) \, dt $$ 
the mean oscillation and 
$$ \| H \|_\infty = \max_{t \in [0,1]} \osc (H_t) $$ 
the maximum oscillation of $H$ on the interval $[0,1]$. 
We call $\| \cdot \|_{(1,\infty)}$ and $\| \cdot \|_\infty$ 
the $L^{(1,\infty)}$-(Hofer) norm and the $L^\infty$-(Hofer) norm 
on the space of smooth time-dependent Hamiltonians respectively (see below for a more precise statement). 
Clearly $\| \cdot \|_{(1,\infty)} \le \| \cdot \|_\infty$, but the converse is false. 
It is easy to see that there are sequences $H_i$ of Hamiltonians such that 
$\| H_i \|_{(1,\infty)} = 1$ for all $i$, but $\| H_i \|_\infty \to \infty$ as $i \to \infty$. 
In fact, one can take any time-independent Hamiltonian $H$ with $\| H \|_\infty = \| H \|_{(1,\infty)} = 1$, 
and define a sequence of reparameterizations (in time) of $H$ so that the constructed sequence has the property above. 

We will consistently use the subscripts (or superscripts) $(1,\infty)$ and $\infty$ to distinguish the two cases, 
and use them to denote any object defined using the one or the other norm. 
We will omit them and write for example $\| \cdot \|$ to denote either one of the two cases. 
That is, when we omit the subscripts (or superscripts), 
the particular statement is true in both the $L^{(1,\infty)}$-case and the $L^\infty$-case. 

For $\phi \in Ham \, (M,\omega)$ we can then define the Hofer norm 
\begin{equation} \label{hamhofer}
\| \phi \| = \inf_{H \mapsto \phi} \| H \| , 
\end{equation}
where $H \mapsto \phi$ means that $H$ generates $\phi$ in the sense that $\phi = \phi_H^1$. 
It is a highly non-trivial fact that this indeed gives a norm on $Ham \, (M,\omega)$, see [1], [8], [2]. 
The above inequality immediately implies $\| \phi \|_{(1,\infty)} \le \| \phi \|_\infty$ 
for $\phi \in Ham \, (M,\omega)$. 
In view of our remark concerning the converse in the case of functions $H \in C^\infty ([0,1] \times M)$, 
the following fact seems rather surprising. 
It is due to Polterovich. 

\begin{lemma}[Lemma 5.1.C, {[9]}] \label{independence}
$\| \phi \|_{(1,\infty)} = \| \phi \|_\infty$ for each $\phi \in Ham \, (M,\omega)$. \qed 
\end{lemma}

Polterovich's proof is constructive. 
Since his arguments will be used in the proof of Theorem \ref{firsttheorem}, 
the proof will be given in section 3. 

Note that both norms $\| \cdot \|$ on $C^\infty ( [0,1] \times M)$ are invariant 
under adding functions that depend only on time $t$. 
We call a Hamiltonian $H$ normalized if 
$$ \int_M H_t \, \omega^n = 0, \hspace{5mm} \mbox{for all } t \in [0,1]. $$ 
Note that normalizing a Hamiltonian in this way does not change $\| \cdot \|$. 
When restricted to the vector space $C_m^\infty ([0,1] \times M)$ 
of normalized smooth time-dependent Hamiltonian functions, 
where $m$ stands for `mean zero', $\| \cdot \|$ indeed defines a norm 
(otherwise it is only a pseudo-norm, since $\| \cdot \|$ vanishes for each function that depends only on time). 
Also note that a normalized Hamiltonian satisfies 
$$ \max_{x \in M} H_t (x) \ge 0, \hspace{5mm} - \min_{x \in M} H_t (x) \ge 0, 
\hspace{5mm} \mbox{for all } t \in [0,1], $$ 
and therefore 
$$ \max_{x \in M} H_t (x) \le \osc (H_t), \hspace{5mm} - \min_{x \in M} H_t (x) \le \osc (H_t), 
\hspace{5mm} \mbox{for all } t \in [0,1]. $$ 
As a consequence, if we denote for all $t \in [0,1]$ 
$$ \| H_t \|_{C^0} = \max_{x \in M} | H_t (x) |, \hspace{1cm} \| H \|_{C^0} = \max_{(t,x)} | H (t,x) |, $$ 
then any normalized Hamiltonian satisfies 
$$ \| H_t \|_{C^0} \le \osc (H_t), \hspace{5mm} \mbox{for all } t \in [0,1], 
\hspace{1cm} \| H \|_{C^0} \le \| H \|_\infty. $$ 

It is easy to see that two Hamiltonians generate the same Hamiltonian path 
if and only if they differ by a function that depends only on time. 
Hence there is a one-one correspondence between smooth Hamiltonian paths 
and generating normalized smooth Hamiltonian functions. 
We therefore frequently consider the Hofer norm as a norm 
on the space $\PP^{ham} (Symp \, (M,\omega), id)$ of smooth Hamiltonian paths. 


\section{The Group of Hamiltonian Homeomorphisms}

Define a metric on the set $\PP^{ham} (Symp \, (M,\omega), id)$ of smooth Hamiltonian paths by 
$$ d_{ham} (\phi_H, \phi_F) 
= \left\| \overline{H} \# F \right\| + \overline{d} (\phi_H, \phi_F) 
= \left\| H - F \right\| + \overline{d} (\phi_H, \phi_F) $$ 
for any two Hamiltonian paths $\phi_H$ and $\phi_F$. 
Here $\overline{H} \# F$ denotes the (normalized) Hamiltonian generating 
the path $(\phi_H)^{-1} \circ \phi_F$, see e.g. [7]. 
Note that we used the above one-one correspondence between smooth Hamiltonian paths 
and smooth normalized Hamiltonian functions. 
On the other hand, 
$$ \overline{d} (\lambda, \mu) = \max_{t \in [0,1]} \overline{d} (\lambda (t), \mu (t) ) $$ 
denotes the $C^0$-distance of any two paths $\lambda$, $\mu$, 
where on the group $Homeo \, (M)$ of homeomorphisms of $M$ $\overline{d}$ is the metric defined by 
$$ \overline{d} (\phi, \psi) = \max \left( d_{C^0} (\phi, \psi), d_{C^0} (\phi^{-1}, \psi^{-1}) \right), $$ 
and where 
$$ d_{C^0} (\phi, \psi) = \max_{x \in M} d \left( \phi (x), \psi (x) \right) $$ 
denotes the standard distance of two maps with respect to a given Riemannian metric $d$ on $M$. 
Since $M$ is compact, the topology induced by $\overline{d}$, which we call the $C^0$-topology, 
agrees with the compact-open topology on $Homeo \, (M)$. 
In particular, this topology is independent of the choice of Riemannian metric $d$. 
Note that we use the notation $\overline{d}$ for both the distance of maps as well as the distance of paths. 
We call the topology induced by the metric $d_{ham}$ the Hamiltonian topology on $\PP^{ham} (Symp \, (M,\omega), id)$, 
or also the $L^{(1,\infty)}$-Hamiltonian topology or $L^\infty$-Hamiltonian topology 
when we need to specify the choice of norm. 
Also see [7] for two different descriptions of the Hamiltonian topology and a detailed discussion of this definition. 

From [7] recall the unfolding map 
\begin{eqnarray*}
(\iota_{ham}, \mbox{Dev}) : \PP^{ham} (Symp \, (M,\omega), id) 
&\to& \PP (Symp \, (M,\omega), id) \times C_m^\infty ([0,1] \times M), \\ 
\hspace{1cm} \lambda = \phi_H &\mapsto& (\lambda, H), 
\end{eqnarray*}
where the topology on the target is the one induced by the product metric of 
the $C^0$-metric $\overline{d}$ on the set $\PP (Symp \, (M,\omega), id)$ of smooth paths 
$\lambda : [0,1] \to Symp \, (M,\omega)$ with $\lambda (0) = id$, 
and the metric induced by $\| \cdot \|$ on $C_m^\infty ([0,1] \times M)$. 
This product embeds into 
$$ \PP (Symp \, (M,\omega), id) \times C_m^\infty ([0,1] \times M) 
\hookrightarrow \PP (Homeo \, (M), id) \times L_m^{(1,\infty)} ([0,1] \times M) $$ 
in the $L^{(1,\infty)}$-case and into 
$$ \PP (Symp \, (M,\omega), id) \times C_m^\infty ([0,1] \times M) 
\hookrightarrow \PP (Homeo \, (M), id) \times C_m^0 ([0,1] \times M) $$ 
in the $L^\infty$-case, 
where $\PP (Homeo \, (M), id)$ denotes the set of continuous paths 
$\lambda : [0,1] \to Homeo \, (M)$ with $\lambda (0) = id$, 
$L_m^{(1,\infty)} ([0,1] \times M)$ denotes the completion of $C_m^\infty ([0,1] \times M)$ 
with respect to the $\| \cdot \|_{(1,\infty)}$-norm, 
and $C_m^0 ([0,1] \times M)$ denotes the vector space of normalized 
(uniformly) continuous functions $H : [0,1] \times M \to \R$. 
The targets of these embeddings are complete metric spaces in the product metrics described above. 

Denote by $\QQ_{(1,\infty)}$ and $\QQ_\infty$ the images of the unfolding maps 
in the $L^{(1,\infty)}$-case and $L^\infty$-case respectively, 
with the induced subspace topologies, which we as well call the Hamiltonian topologies 
(see [Remark 3.17 (2) in 7] for an explanation of this terminology). 
Denote by $\overline{\QQ}_{(1,\infty)}$ and $\overline{\QQ}_\infty$ their metric completions respectively. 
That is, a pair $(\lambda, H)$, where $\lambda : [0,1] \to Homeo \, (M)$ is a continuous path with $\lambda (0) = id$ 
and $H : [0,1] \times M \to \R$ is an $L^{(1,\infty)}$-function (in the $L^{(1,\infty)}$-case) 
or a continuous function (in the $L^\infty$-case), 
lies in $\overline{\QQ}_{(1,\infty)}$ or $\overline{\QQ}_\infty$ respectively, 
if and only if there is a sequence $(\phi_{H_i}, H_i)$, 
where the $H_i$ are smooth normalized Hamiltonians and the $\phi_{H_i}$ are the corresponding smooth Hamiltonian paths, 
such that the sequence $H_i$ converges to $H$ in the Hofer norm $\| \cdot \|$ 
and the sequence $\phi_{H_i}$ converges to $\lambda$ in the $C^0$-metric. 

Note that given a sequence $(\phi_{H_i}, H_i)$ as above, any subsequence has the same limit $(\lambda, H)$. 
In particular, given any (decreasing) sequence $\e_i > 0$ of positive numbers with $\e_i \to 0$ as $i \to \infty$, 
by passing to a subsequence we may assume that the given sequence satisfies $\| H_j - H_k \| < \e_i$ 
and $\overline{d} (\phi_{H_j}, \phi_{H_k}) < \e_i$ for all $j, k \ge i$, for all $i$. 
Similarly, we may assume that the given sequence satisfies $\| H_j - H \| < \e_i$ 
and $\overline{d} (\phi_{H_j}, \lambda) < \e_i$ for all $j \ge i$ and all $i$, 
or any combination of the two, e.g. we may assume that the given sequence satisfies 
$\| H_j - H_k \| < \e_i$ for all $j, k \ge i$ 
and $\overline{d} (\phi_{H_j}, \lambda) < \e_i$ for all $ j \ge i$, and for all $i$. 
It is often convenient to consider a sequence $\e_i$ such that $\sum \e_i$ converges, 
or $\sum_{i=j}^k \e_i \to 0$ as $k \ge j \to \infty$, 
for example $\e_i = \frac{1}{2^i}$. 

Recall that for any sequence given as above, each $H_i$ is normalized. 
In particular (see section 1), unlike in the $L^{(1,\infty)}$-case, 
in the $L^\infty$-case $H_i$ converges uniformly on $[0,1] \times M$ to a uniformly continuous function $H$ 
(and therefore the metric space considered above is indeed complete). 
This gives an obvious advantage for many arguments in the $L^\infty$-case. 
The Reparameterization Lemma [Lemma 3.21 in 7], 
the Structure Theorem for topological Hamiltonians [4, or Theorem 2.6 in 6], 
and similar arguments in the $L^{(1,\infty)}$-case 
were specifically designed to make up for the lack of continuity. 

The evaluation map 
$$ ev_1 : \PP^{ham} (Symp \, (M,\omega), id) \to Ham \, (M,\omega), \hspace{1cm} \lambda \mapsto \lambda(1), $$ 
induces a well-defined evaluation map 
\begin{equation} \label{oneinftyevalmap}
\overline{ev}_1^{\QQ_{(1,\infty)}} : \overline{\QQ}_{(1,\infty)} \to Homeo \, (M), 
\hspace{1cm} \lambda \mapsto \lambda(1), 
\end{equation}
see [7]. 
The same arguments apply to the $L^\infty$-case, so we also get an evaluation map 
\begin{equation} \label{inftyevalmap}
\overline{ev}_1^{\QQ_\infty} : \overline{\QQ}_\infty \to Homeo \, (M), \hspace{1cm} \lambda \mapsto \lambda(1). 
\end{equation}

\begin{definition}[Definition 3.14, {[7]}] \label{definitionhameo}
The sets $Hameo_{(1,\infty)} \, (M,\omega)$ and $Hameo_\infty \, (M,\omega)$ 
of ($L^{(1,\infty)}$- and $L^\infty$-)Hamiltonian homeomorphisms are defined as the images 
$$ Hameo \, (M,\omega) = \, \mbox{Im} \, \left( \overline{ev}_1^{\QQ} \right) \subset Homeo \, (M) $$ 
of the above evaluation maps (\ref{oneinftyevalmap}) and (\ref{inftyevalmap}) respectively, 
equipped with the subspace topology from $Homeo \, (M)$, i.e., the $C^0$-topology. 
By definition the ($L^{(1,\infty)}$- and $L^\infty$-)Hamiltonian topology is the strongest topology on these sets 
such that the evaluation maps (\ref{oneinftyevalmap}) and (\ref{inftyevalmap}) are continuous respectively. 
The resulting topological spaces are denoted by 
$\HH ameo_{(1,\infty)} \, (M,\omega)$ and $\HH ameo_\infty \, (M,\omega)$ respectively. 
\end{definition}

In other words, $h \in Hameo \, (M,\omega)$ if and only if 
there exists a Cauchy sequence $(\phi_{H_i}, H_i) \in \QQ$ in the Hamiltonian topology 
such that $\phi_{H_i}^1 \to h$ in the $C^0$-metric. 
In particular, $Ham \, (M,\omega) \subset Hameo \, (M,\omega)$, 
where we recall the convention that if we drop the subscripts 
the statement holds in both the $L^{(1,\infty)}$-case and the $L^\infty$-case. 

\begin{theorem}[Corollary 3.25, {[7]}]
$Hameo \, (M,\omega) \subset Homeo \, (M)$ is a topological subgroup. \qed 
\end{theorem}

In the $L^{(1,\infty)}$-case the proof is given in [7]. 
The proof in the $L^\infty$-case is essentially the same. 
In fact, some arguments even simplify since, as remarked above, Cauchy sequences in the $L^\infty$-Hofer norm 
converge uniformly and therefore have uniformly continuous limits. 

Denote by $Sympeo \, (M,\omega) = \overline{Symp \, (M,\omega)} \subset Homeo \, (M)$ 
the group of symplectic homeomorphisms of $(M,\omega)$, 
where the closure is with respect to the $C^0$-topology. 
See [7] for a foundational study of $Sympeo \, (M,\omega)$. 
The following theorems are proved in [7] in the $L^{(1,\infty)}$-case. 
Again, their proofs can easily be adapted to the $L^\infty$-case. 

\begin{theorem}[Theorem 4.4, {[7]}] \label{normalsubgp}
$Hameo \, (M,\omega)$ is a normal subgroup of $Sympeo \, (M,\omega)$. \qed 
\end{theorem}

\begin{theorem}[Theorem 4.5, {[7]}]
$Hameo \, (M,\omega)$ and $\HH ameo \, (M,\omega)$ are path-connected, 
and $\HH ameo \, (M,\omega)$ is locally path-connected. \qed 
\end{theorem}

The question of local path-connectedness of $Hameo \, (M,\omega)$ is open. 

Recall from section 1 that $\| \cdot \|_{(1,\infty)} \le \| \cdot \|_\infty$. 
Therefore any Cauchy sequence in the $L^\infty$-Hamiltonian topology 
is also a Cauchy sequence in the $L^{(1,\infty)}$-Hamiltonian topology. 
In particular, $Hameo_\infty \, (M,\omega) \subset Hameo_{(1,\infty)} \, (M,\omega)$. 
In view of Theorem \ref{normalsubgp} we in fact have 

\begin{proposition} \label{subgroup}
$Hameo_\infty \, (M,\omega)$ forms a normal subgroup of $Hameo_{(1,\infty)} \, (M,\omega)$. \qed 
\end{proposition}

Polterovich's Lemma \ref{independence} motivates the question whether equality holds in Proposition \ref{subgroup}. 
Note that if the inclusion were proper, 
this would negatively answer the simpleness question of $Hameo_{(1,\infty)} \, (M,\omega)$, 
which was raised in [7], 
and would also give a (simple) proof of nonsimpleness of the kernel of the mass flow homomorphism 
on a smooth closed orientable surface, see [section 5 in 7]. 
As one therefore expects, equality indeed holds. 
This will be proved in the next section. 

The projections $\overline{\iota}_{ham}^{\QQ}$ of $\overline{\QQ}_{(1,\infty)}$ 
and $\overline{\QQ}_\infty$ onto the first factor 
yield the spaces of ($L^{(1,\infty)}$- and $L^\infty$-)topological Hamiltonian paths, 
denoted by $\PP^{ham}_{(1,\infty)} (Sympeo \, (M,\omega), id)$ 
and $\PP^{ham}_\infty (Sympeo \, (M,\omega), id)$ respectively. 
And the projections $\overline{\Dev^{\QQ}}$ onto the second factor 
yield the spaces of (normalized) ($L^{(1,\infty)}$- and $L^\infty$-)topological Hamiltonian functions, 
denoted by $\HH_m^{(1,\infty)} ([0,1] \times M)$ and $\HH_m^\infty ([0,1] \times M)$ respectively. 
It follows from the following statements that, at least in the $L^\infty$-case, 
the above one-one correspondence between smooth Hamiltonian paths 
and normalized smooth Hamiltonian functions extends to the topological Hamiltonian category. 

\begin{theorem}[Uniqueness Theorem, {[10]}] \label{functionuniqueness}
The map $\overline{\iota}_{ham}^{\QQ_\infty} : 
\overline{\QQ_\infty} \to \PP^{ham}_\infty (Sympeo \, (M,\omega), id)$ is one-one. 
In other words, if $(\lambda, H)$ and $(\lambda, H') \in \overline{\QQ_\infty}$, then $H = H'$. \qed 
\end{theorem}

Note that the normalization condition is crucial in the above statement. 

\begin{theorem}[Theorem 3.1, {[6]}] \label{pathuniqueness}
The map $\overline{\Dev^{\QQ}} : \overline{\QQ} \to \HH_m ([0,1] \times M)$ is one-one. 
In other words, if $(\lambda, H)$ and $(\lambda', H) \in \overline{\QQ}$, then $\lambda = \lambda'$. \qed 
\end{theorem}

Theorem \ref{functionuniqueness} was proved in [10], and later more details were added in the unpublished preprint [4], 
which attempts to generalize the result in [10] to the $L^{(1,\infty)}$-case. 
A local version of Theorem \ref{functionuniqueness} can be found in [5]. 
The proof of Theorem \ref{pathuniqueness} is given in [6] in the $L^{(1,\infty)}$-case, 
but applies to the $L^\infty$-case without changes. 
This one-one correspondence allows to extend Hamiltonian dynamics to the topological Hamiltonian category, see [6]. 

Note that by the argument before Proposition \ref{subgroup}, we have 

\begin{proposition} \label{inclusion}
$\HH_m^\infty ([0,1] \times M) \subset \HH_m^{(1,\infty)} ([0,1] \times M)$ and \\ 
$\PP^{ham}_\infty (Sympeo \, (M,\omega), id) \subset \PP^{ham}_{(1,\infty)} (Sympeo \, (M,\omega), id)$. \qed 
\end{proposition}

We finish this section expanding the discussion from [6] on the length or Hofer norm of a topological Hamiltonian path 
and the Hofer norm of a Hamiltonian homeomorphism. 
We will assume a uniqueness result as in Theorem \ref{functionuniqueness}. 
That is, the following applies to the $L^\infty$-case, 
but will also apply to the $L^{(1,\infty)}$-case if the analog to Theorem \ref{functionuniqueness} can be established. 

Let $\lambda \in \PP^{ham} (Sympeo \, (M,\omega), id)$ be a topological Hamiltonian path. 
By definition, there is a sequence $(\phi_{H_i}, H_i) \in \QQ$ so that $\phi_{H_i} \to \lambda$ in the $C^0$-metric 
and $H_i$ converges in Hofer's norm to a topological Hamiltonian $H$. 
We then define the length or Hofer norm of the topological Hamiltonian path $\lambda$ by 
\begin{equation} \label{leng}
\leng (\lambda) = \lim_{i \to \infty} \leng (\phi_{H_i}) = \lim_{i \to \infty} \| H_i \| = \| H \| . 
\end{equation}
The Uniqueness Theorem \ref{functionuniqueness} 
(and an application of the triangle inequality) implies that this definition is well-defined. 
In particular, this agrees with the usual definition $\leng (\phi_H) = \| H \|$ 
of the length of a smooth Hamiltonian path. 
For $h \in Hameo \, (M,\omega)$ a Hamiltonian homeomorphism, we define the Hofer norm of $h$ to be 
\begin{equation} \label{hameohofer}
\| h \| = \inf \left \{ \leng (\lambda) \mid \lambda \in \PP^{ham} (Sympeo \, (M,\omega), id), 
\; \overline{ev}_1^{\QQ} (\lambda) = h \right \} . 
\end{equation}
This is obviously well-defined. 
It is not too hard to see that the Hofer norm $\| \cdot \|$ is continuous 
with respect to the Hamiltonian topology (Definition \ref{definitionhameo}) on $Hameo \, (M,\omega)$. 
Moreover, it has the same basic properties as the usual Hofer norm on $Ham \, (M,\omega)$, i.e., 
it is symmetric, bi-invariant, satisfies the triangle inequality and symplectic invariance, 
and it is nondegenerate (and therefore indeed defines a norm) [6]. 

Let us denote by $\| \cdot \|_{Ham}$ the usual Hofer norm (\ref{hamhofer}) 
on $Ham \, (M,\omega)$ defined in the beginning 
and by $\| \cdot \|_{Hameo}$ the Hofer norm (\ref{hameohofer}) on $Hameo \, (M,\omega)$ defined above. 
Then both norms are defined on $Ham \, (M,\omega) \subset Hameo \, (M,\omega)$. 
To avoid confusion we will always use these subscripts to denote 
either of the Hofer norms of a Hamiltonian diffeomorphism, 
but omit the subscript \textit{Hameo} when we mean the Hofer norm of a Hamiltonian homeomorphism. 
We clearly have $\| \cdot \|_{Hameo} \le \| \cdot \|_{Ham}$. 
It seems likely that they are in fact equal. 

\begin{question} \label{question}
Is $\| \cdot \|_{Hameo} = \| \cdot \|_{Ham}$ on $Ham \, (M,\omega)$? 
\end{question}

The answer to this question is not known. 
The difficulty is that if, for some $\phi \in Ham \, (M,\omega)$, 
$\lambda$ is a topological Hamiltonian path whose length approximates $\| \phi \|_{Hameo}$, 
then $\lambda$ can be approximated by smooth Hamiltonian paths in the Hamiltonian topology. 
However, the end points of these paths are in general different from $\phi$, 
and therefore these paths are not admissable to compute the Hofer norm $\| \phi \|_{Ham}$ of $\phi$. 
Note that there is a `short' topological Hamiltonian path from each such end point to $\phi$, 
but there need not be such a smooth path. 
So a related question is whether we can choose above sequence of smooth Hamiltonian paths 
so that the end point equals $\phi$ for each path. 
This problem seems to lie at the heart of topological Hamiltonian geometry. 

In view of Theorem \ref{firsttheorem} and Lemma \ref{independence}, 
one can then ask whether we have $\| \cdot \|_\infty = \| \cdot \|_{(1,\infty)}$ 
for the Hofer norms (\ref{hameohofer}) on $Hameo \, (M,\omega)$. 
Note that if the answer to Question \ref{question} were affirmative, 
then this together with Lemma \ref{independence} 
and the continuity of the Hofer norm with respect to the Hamiltonian topology 
would immediately imply this equality. 
It is not too hard to show that (still assuming uniqueness of topological Hamiltonians) equality indeed holds 
by proving an analog to Lemma \ref{technicallemma} in the next section for topological Hamiltonians. 

The above would in particular imply that given $h \in Hameo \, (M,\omega)$ and $\e > 0$, 
there exists a path $\lambda$ which is a topological Hamiltonian path 
in both the $L^{(1,\infty)}$-sense and the $L^\infty$-sense and such that 
$$ \| h \|_\infty = \| h \|_{(1,\infty)} \le \lengoneinf (\lambda) \le \lenginf (\lambda) < \| h \|_\infty + \e . $$ 
In other words, Theorem \ref{firsttheorem} states that each end point of an $L^{(1,\infty)}$-topological Hamiltonian path 
is also the end point of some (possibly different) $L^\infty$-topological Hamiltonian path, 
and if uniqueness holds in both cases (so that the Hofer norm is well-defined), 
then the same statement holds for `short' topological Hamiltonian paths as well. 


\section{Proof of the Main Theorem}

In this section we prove the main theorem, which we restate for convenience 

\begin{theorem} \label{maintheorem}
$Hameo_\infty \, (M,\omega) = Hameo_{(1,\infty)} \, (M,\omega)$. \qed 
\end{theorem}

By Proposition \ref{subgroup} we already have the inclusion 
$Hameo_\infty \, (M,\omega) \subset Hameo_{(1,\infty)} \, (M,\omega)$. 
The converse is more delicate. 
To deduce it, we will first prove a series of lemmas. 
Our proof is mainly based on Polterovich's Lemma \ref{independence}. 
In fact, Polterovich proved the slightly stronger Lemma \ref{technicallemma} below. 
Lemma \ref{independence} will be an immediate consequence of this more technical result. 

In the following we want to allow more general smooth Hamiltonian paths that are not necessarily based at the identity. 
That is, unless explicit mention is made to the contrary, 
we consider paths $\lambda = \phi \circ \phi_H$, 
where $\phi_H$ is a smooth Hamiltonian path in the previous sense with $\phi_H^0 = id$, and $\phi \in Ham \, (M,\omega)$. 
It is easy to see that $\lambda$ solves Hamilton's equation 
with Hamiltonian $H \circ \phi$ and initial condition $\lambda (0) = \phi$. 
We therefore call the Hamiltonian $K = H \circ \phi$ the generating Hamiltonian of $\lambda$. 
Note that we could also work with paths of the form $\lambda = \phi_H \circ \phi$, 
which solve Hamilton's equation with Hamiltonian $H$ and initial condition $\lambda (0) = \phi$. 
It turns out that the former will be more convenient for the computations below. 
We will often simply write $\lambda = \phi_K$. 
If one path starts where another one ends, we may consider their concatenation, 
and if both paths are boundary flat (see below for the definition), then that concatenation is a smooth Hamiltonian path 
whose generating Hamiltonian agrees with (a reparameterization of) the Hamiltonian of the first path for some time 
and with (a reparameterization of) the Hamiltonian of the second path for the remaining time. 
Of course, when computing the Hofer norm of an element $\phi \in Ham \, (M,\omega)$, 
we only allow paths $\phi_H$ with $\phi_H^0 = id$ and $\phi_H^1 = \phi$. 

\begin{lemma} \label{technicallemma}
Let $H : [0,1] \times M \to \R$ be a smooth normalized Hamiltonian function 
generating the smooth Hamiltonian path $\phi_H : t \in [0,1] \mapsto \phi_H^t$. 
Let $\e > 0$ be given. 
Then there exists a smooth normalized Hamiltonian function $F : [0,1] \times M \to \R$ such that the following holds 
\begin{description}
\item [(i)] $\phi_F^0 = \phi_H^0$ and $\phi_F^1 = \phi_H^1$, 
\item [(ii)] $\| F \|_\infty < \| H \|_{(1,\infty)} + \e$, and 
\item [(iii)] $\overline{d} (\phi_F, \phi_H^0) < \overline{d} (\phi_H, \phi_H^0) + \e$. 
\end{description}
In (iii) $\phi_H^0$ denotes the constant path $t \mapsto \phi_H^0$. \qed 
\end{lemma}

Recall the following reparameterization procedure (see for example [7] or [9]). 
Let $H : [0,1] \times M \to \R$ be a smooth Hamiltonian function 
generating the smooth Hamiltonian path $\phi_H : t \mapsto \phi_H^t$. 
For any smooth function $\zeta : [0,1] \to [0,1]$ the reparameterized path $\phi_H^\zeta : t \mapsto \phi_H^{\zeta (t)}$ 
is generated by the smooth Hamiltonian function $H^\zeta$ 
given by the formula $H^\zeta (t,x) = \zeta' (t) H ( \zeta (t), x)$. 
If $\zeta (0) = 0$ and $\zeta (1) = 1$, then the time-zero maps and time-one maps coincide respectively, 
that is, $\phi_{H^\zeta}^0 = \phi_H^{\zeta (0)} = \phi_H^0$ and $\phi_{H^\zeta}^1 = \phi_H^{\zeta (1)} = \phi_H^1$. 
Moreover, if $H$ is normalized then so is $H^\zeta$. 
We refer to the function $\zeta$ as the reparameterization function. 

\begin{proof}[Proof of Lemma \ref{technicallemma}]
We first consider the path $t \mapsto \phi_K^t = \psi^t \circ \phi_H^t$, 
where $t \mapsto \psi^t$ is a loop in $Ham \, (M,\omega)$, $\psi^0 = \psi^1 = id$. 
Clearly $\phi_K^0 = \phi_H^0$ and $\phi_K^1 = \phi_H^1$. 
We may choose the loop $\psi^t$ such that it is arbitrarily close to the constant loop $id$ in the $C^0$-metric, 
its generating Hamiltonian is arbitrarily small in the $\| \cdot \|_{(1,\infty)}$-norm, 
and such that $\osc (K_t) \not= 0$ for all $t \in [0,1]$, see [section 5.2 in 9]. 
Therefore we may choose the Hamiltonian $K$ such that 
$$ \| K \|_{(1,\infty)} < \| H \|_{(1,\infty)} + \frac{\e}{2}, 
\hspace{1cm} \overline{d} (\phi_K, \phi_H^0) < \overline{d} (\phi_H, \phi_H^0) + \e . $$ 
To see the second inequality, write 
$$ \overline{d} (\phi_K, \phi_H^0) \le \overline{d} (\phi_K, \phi_H) + \overline{d} (\phi_H, \phi_H^0), $$ 
and note that the first term on the right-hand side of the inequality can be made as small as we want 
since the set of continuous paths $\PP (Homeo \, (M), id)$ forms a topological group in the $C^0$-topology [7]. 

We may normalize $K$ if necessary without losing any of the above properties. 
Now reparameterize $K$ to $K^\zeta$, where $\zeta$ is the inverse of 
(here we use $\osc (K_t) \not= 0$ for all $t \in [0,1]$) 
$$ \xi : [0,1] \to [0,1], \hspace{1cm} t \mapsto \frac{\int_0^t \osc (K_s) \, ds}{\int_0^1 \osc (K_s) \, ds} . $$ 
Note that $\zeta$ fixes $0$ and $1$, so that $\phi_{K^\zeta}$ has the same end points as $\phi_K$. 
Then 
$$ \zeta' (t) = \frac{\int_0^1 \osc (K_s) \, ds}{\osc K_{\zeta (t)}} . $$ 
Hence for every $t$ 
$$ \osc \left( (K^\zeta)_t \right) = \zeta' (t) \, \osc (K_{\zeta (t)}) 
= \int_0^1 \osc (K_s) \, ds = \| K \|_{(1,\infty)} , $$ 
and therefore 
$$ \left\| K^\zeta \right\|_\infty = \| K \|_{(1,\infty)} < \| H \|_{(1,\infty)} + \frac{\e}{2}. $$ 
Now $\zeta$ (and therefore $K^\zeta$ and $\phi_{K^\zeta}$) may not be smooth but only $C^1$. 
So we approximate $\zeta$ in the $C^1$-topology by a smooth diffeomorphism $\rho$ of $[0,1]$ that also fixes $0$ and $1$ 
to obtain a smooth normalized Hamiltonian $F = K^\rho$ with $\| F \|_{\infty} < \| K^\zeta \|_\infty + \frac{\e}{2}$. 
Then $F$ clearly satisfies (i) and (ii). 
Since $\phi_F$ is just a reparameterization of $\phi_K$ we also have 
$$ \overline{d} (\phi_F, \phi_H^0) = \overline{d} (\phi_K, \phi_H^0) < \overline{d} (\phi_H, \phi_H^0) + \e . $$ 
That proves (iii) and hence finishes the proof. 
\end{proof}

\begin{proof}[Proof of Lemma \ref{independence}]
For every Hamiltonian $H$ we have $\| H \|_{(1,\infty)} \le \| H \|_\infty$. 
So the inequality $\| \phi \|_{(1,\infty)} \le \| \phi \|_\infty$ is obvious. 
For the converse, let $\e > 0$ be arbitrary. 
Choose a Hamiltonian $H$ generating $\phi$ such that $\| H \|_{(1,\infty)} \le \| \phi \|_{(1,\infty)} + \e$. 
By Lemma 3.2, we can find a Hamiltonian $F$ generating $\phi$ such that 
$\| F \|_\infty < \| H \|_{(1,\infty)} + \e \le \| \phi \|_{(1,\infty)} + 2 \e$. 
But then $\| \phi \|_\infty \le \| F \|_\infty < \| \phi \|_{(1,\infty)} + 2 \e$. 
Since $\e > 0$ was arbitrary, this implies $\| \phi \|_\infty \le \| \phi \|_{(1,\infty)}$. 
This completes the proof. 
\end{proof}

The next lemma is proved in [3] and in the form stated here in [7]. 

\begin{lemma} [$L^{(1,\infty)}$-Approximation Lemma, Lemma 5.2, {[3]}, Lemma A.4, {[7]}] \label{oneinftyapproximation}
Let $H : [0,1] \times M \to \R$ be a smooth normalized Hamiltonian function 
generating the smooth Hamiltonian path $\phi_H : t \mapsto \phi_H^t$. 
Then given any $\e > 0$, there exists a smooth normalized Hamiltonian function $F : [0,1] \times M \to \R$ such that 
\begin{itemize}
\item $F$ (and hence $\phi_F$) is boundary flat, that is, 
there exists $\delta > 0$ such that $F_t \equiv 0$ for $0 \le t \le \delta$ and $1 - \delta \le t \le 1$, 
\item $\phi_F^0 = \phi_H^0$ and $\phi_F^1 = \phi_H^1$, 
\item $\| F - H \|_{(1,\infty)} < \e$, and 
\item $\overline{d} (\phi_F, \phi_H) < \e$. \qed 
\end{itemize}
\end{lemma}

In other words, any smooth Hamiltonian path can be approximated 
by a boundary flat smooth Hamiltonian path in the $L^{(1,\infty)}$-Hamiltonian topology. 
It is easy to see that this approximation procedure fails in the $L^\infty$-norm 
unless (after normalizing) $H_0 \equiv H_1 \equiv 0$, 
since for any boundary flat Hamiltonian F we have 
$\| F - H \|_\infty \ge \max \left( \| F_0 - H_0 \|_{C^0}, \| F_1 - H_1 \|_{C^0} \right) 
= \max ( \| H_0 \|_{C^0}, \| H_1 \|_{C^0})$. 
Since for a given Cauchy sequence we cannot expect $H_{i,0}$ and $H_{i,1} \to 0$ as $i \to \infty$ in general, 
there can be no boundary flattening procedure for Cauchy sequences in the $L^\infty$-case 
as for Cauchy sequences in the $L^{(1,\infty)}$-case. 
The continuity of this boundary flattening procedure is one of the main advantages of the $L^{(1,\infty)}$-norm. 
The boundary flattening procedure is important for 
extending Yong-Geun Oh's spectral invariants to the topological Hamiltonian category, see [6], 
and is the main reason we adopt the $L^{(1,\infty)}$-norm in general. 
We would also like to remark that the $L^{(1,\infty)}$-norm seems more natural 
in the context of Floer theory and the theory of currents, see [7], [4], [6], 
and for the notion of the length of a (topological) Hamiltonian path (section 2). 

However, there is a version of this reparameterization lemma in the $L^\infty$-case 
that can be applied to sequences of Hamiltonians converging to zero. 
Its proof is along the same lines as the proof of the $L^{(1,\infty)}$-Approximation Lemma. 

\begin{lemma} \label{inftyapproximation}
Let $H : [0,1] \times M \to \R$ be a smooth normalized Hamiltonian function 
generating the smooth Hamiltonian path $\phi_H : t \mapsto \phi_H^t$. 
There exists a positive constant $C$ that depends only on $H$ such that, 
given any $\e > 0$, 
there exists a smooth normalized Hamiltonian function $F : [0,1] \times M \to \R$ such that 
\begin{itemize}
\item $F$ (and hence $\phi_F$) is boundary flat, 
\item $\phi_F^0 = \phi_H^0$ and $\phi_F^1 = \phi_H^1$, 
\item $\| F - H \|_\infty \le 2 \| H \|_\infty + C \e $, 
and in particular, $\| F \|_\infty \le 3 \| H \|_\infty + C \e$, 
\item $\overline{d} (\phi_F, \phi_H) < \e$. \qed 
\end{itemize}
\end{lemma}

\begin{proof}
We choose a smooth reparameterization function $\zeta : [0,1] \to [0,1]$ with the following properties 
\begin{itemize}
\item $\zeta \equiv 0$ near $t = 0$ and $\zeta \equiv 1$ near $t = 1$, 
\item $\| \zeta - id \|_{C^0} <\e$, and 
\item $0 \le \zeta' (t) \le 2$ for all $t \in [0,1]$. 
\end{itemize}
Let $F = H^\zeta$. 
Then $F$ satisfies the first two properties. 
Since $H$ is smooth and $M$ is compact there exists a constant $L < \infty$ 
such that $\| H_t - H_s \|_{C^0} < L | t - s |$ for all $t, s \in [0,1]$. 
Since both $H$ and $F$ are normalized we have for each $t \in [0,1]$ 
\begin{eqnarray*}
0 \le \max_{x \in M} \left( ( H^\zeta )_t - H_t \right) 
&\le& \zeta' (t) \left\| H_{\zeta (t)} - H_t \right\|_{C^0} + \left| \zeta'(t) - 1 \right| \| H_t \|_{C^0} \\ 
&\le& 2 L \| \zeta - id \|_{C^0} + \osc (H_t) \\ 
&<& 2 L \e + \| H \|_\infty . 
\end{eqnarray*}
The same estimate holds for $- \min_{x\in M} \left( ( H^\zeta )_t - H_t \right)$. 
Therefore, $\| F - H \|_\infty < 4 L \e + 2 \| H \|_\infty$. 
That shows that the third property holds. 
The last statement is proved in the same way as in the $L^{(1,\infty)}$-Approximation Lemma. 
\end{proof}

Before completing the proof of Theorem 3.1, we wish to consider another way of reparameterizing a Hamiltonian $H$. 
Given $0 \le a < b \le 1$, and a smooth Hamiltonian function $H$ defined on $[0,1] \times M$, 
we denote by $\zeta_{a,b} : [a,b] \to [0,1]$ the unique linear function with $\zeta (a) = 0$ and $\zeta (b) = 1$, 
and by $H^{\zeta_{a,b}}$ the reparameterized smooth Hamiltonian defined on $[a,b] \times M$. 
Of course if $H$ is normalized then so is $H^{\zeta_{a,b}}$, 
and if $H$ is boundary flat then again so is $H^{\zeta_{a,b}}$. 
Obviously, $\| H^{\zeta_{a,b}} \|_{(1,\infty)} = \| H \|_{(1,\infty)}$ 
and $\| H^{\zeta_{a,b}} \|_\infty = \frac{1}{b-a} \| H \|_\infty$. 

It only remains to show the inclusion $Hameo_{(1,\infty)} \, (M,\omega) \subset Hameo_\infty \, (M,\omega)$. 
Let $h \in Hameo_{(1,\infty)} \, (M,\omega)$. 
By definition, there exists a sequence $(\phi_{H_i}, H_i)$ of smooth normalized Hamiltonian functions $H_i$ 
generating the smooth Hamiltonian paths $\phi_{H_i}$ such that 
$(\phi_{H_i}, H_i)$ converges in the $L^{(1,\infty)}$-Hamiltonian topology and $\phi_{H_i} \to h$ in the $C^0$-metric. 
As remarked in section 1, we cannot expect the sequence $H_i$ to be Cauchy in the $L^\infty$-norm in general. 
Our goal is to modify the sequence $(\phi_{H_i}, H_i)$ to a sequence 
which is Cauchy in the $L^\infty$-Hamiltonian topology. 
Our strategy will be as follows. 
The given sequence gives a `short' path from the end point $\phi_i$ of the path $\phi_{H_i}$ 
to the end point $\phi_{i+1}$ of the path $\phi_{H_{i+1}}$ for all $i$. 
We will construct a sequence $\phi_{F_i}$ of Hamiltonian paths so that $\phi_{F_{i+1}}$ coincides with 
its predecessor $\phi_{F_i}$ for some time, followed by the path from $\phi_i$ to $\phi_{i+1}$. 
We will have to apply Lemma \ref{technicallemma} to pass from the $L^{(1,\infty)}$-norm to the $L^\infty$-norm, 
and we have to make the pieces that we paste together boundary flat 
so that the elements of the constructed sequence are smooth. 
Along the way, we will have to keep track of the closeness of the paths 
and of their Hamiltonians in the Hamiltonian topology. 
Note that the image in $Homeo \, (M)$ of the limit path of the modified sequence 
is very different from the image of the limit path of the original sequence. 
We cannot apply Lemma \ref{technicallemma} directly to the sequence $H_i$ 
and expect these Hamiltonians and the paths they generate to be Cauchy in the $L^\infty$-Hamiltonian topology in general. 

\begin{proof}[Proof of Theorem \ref{maintheorem}]
As remarked before it only remains to show 
the inclusion $Hameo_{(1,\infty)} \, (M,\omega) \subset Hameo_\infty \, (M,\omega)$. 
Let $h \in Hameo_{(1,\infty)} \, (M,\omega)$. 
By definition, there exists a sequence $(\phi_{H_i}, H_i)$ of smooth normalized Hamiltonian functions $H_i$ 
generating the smooth Hamiltonian paths $\phi_{H_i}$ (with $\phi_{H_i}^0 = id$) such that 
\begin{itemize}
\item $\overline{d} (\phi_{H_i}, \phi_{H_j}) \to 0$, as $i, j \to \infty$, 
\item $\left\| \overline{H_i} \# H_j \right\|_{(1,\infty)} 
= \| H_j - H_i \|_{(1,\infty)} \to 0$, as $i, j \to \infty$, and 
\item $\overline{d} (\phi_{H_i}^1, h) \to 0$ as $i \to \infty$. 
\end{itemize}
Denote by $\lambda$ the $C^0$-limit of the sequence of paths $\phi_{H_i}$. 

Let $\e_i > 0$ be a decreasing sequence of real numbers such that $\e_i \to 0$ as $i \to \infty$. 
Since $h$ is uniformly continuous, there exists a sequence $\delta_i > 0$ such that for all $i$, 
$d (h (x), h (y) ) < \e_i$ for all $x, y \in M$ with $d (x, y) < \delta_i$. 
W.l.o.g. we may assume that $\delta_i \le \e_i$ for all $i$. 
As explained in section 2, by passing to a subsequence if necessary we may assume that 
$$ \left\| \overline{H_i} \# H_{i+1} \right\|_{(1,\infty)} < \delta_i \le \e_i, \hspace{5mm} 
\overline{d} (\phi_{H_i}, \lambda) < \delta_i, \hspace{5mm} \mbox{for all } i . $$ 
Assume this is done. 
We will specify the sequence $\e_i$ later in the proof. 

For convenience denote by $H_0$ the Hamiltonian $H_0 = 0$, which generates the constant loop $id$. 
Define the sequence of smooth Hamiltonians $K_i$ by 
$K_i = \left( \overline{H_{i-1}} \# H_i \right) \circ \phi_{i-1}$ for all $i > 1$. 
Since the Hamiltonians $H_i$ are normalized for all $i$, 
the $K_i$ are normalized as well, and the Hamiltonian paths they generate can be chosen to be the paths 
$\phi_{K_i} = \phi_{i-1} \circ \left( \phi_{H_{i-1}} \right)^{-1} \circ \phi_{H_i}$ 
from $\phi_{i-1}$ to $\phi_i$ for all $i > 1$. 
Here and in the following we denote by $\phi_i = \phi_{H_i}^1$ the diffeomorphism itself 
or the constant path $t \mapsto \phi_{H_i}^1$, 
and similarly for the identity $id$. 

By assumption, $\| K_i \|_{(1,\infty)} < \e_{i-1}$ for all $i > 1$. 
Moreover, we claim that the assumption on the sequence $\phi_{H_i}$ implies 
$$ \overline{d} ( \phi_{K_i}, \phi_{i-1} ) \le 4 \e_{i-1} $$ 
for all $i > 1$. 
Again recall that the set of continuous paths $\PP (Homeo (M), id)$ 
forms a topological group in the $C^0$-topology. 
We will use a similar argument as in the proof of this fact to see the above inequality. 
Namely, note that by definition 
\begin{eqnarray*}
\overline{d} ( \phi_{K_i}, \phi_{i-1} ) 
&=& \max \left( d_{C^0} \left( \phi_{K_i}, \phi_{i-1} \right), 
d_{C^0} \left( (\phi_{K_i})^{-1}, \phi_{i-1}^{-1} \right) \right) \\ 
&=& \max \left( d_{C^0} \left( \phi_{i-1} \circ \left( \phi_{H_{i-1}} \right)^{-1} \circ \phi_{H_i}, \phi_{i-1} \right), 
d_{C^0} \left( \left( \phi_{H_i} \right)^{-1} \circ \phi_{H_{i-1}} \circ \phi_{i-1}^{-1}, 
\phi_{i-1}^{-1} \right) \right). 
\end{eqnarray*}
For the second term, we use that the metric $d_{C^0}$ is right-invariant to see that 
\begin{eqnarray*}
d_{C^0} \left( \left( \phi_{H_i} \right)^{-1} \circ \phi_{H_{i-1}} \circ \phi_{i-1}^{-1} , \phi_{i-1}^{-1} \right) 
&=& d_{C^0} \left( \left( \phi_{H_i} \right)^{-1} \circ \phi_{H_{i-1}} , id \right) \\ 
&=& d_{C^0} \left( \left( \phi_{H_i} \right)^{-1}, \left( \phi_{H_{i-1}} \right)^{-1} \right) \\ 
&\le& \overline{d} \left( \left( \phi_{H_i} \right)^{-1}, \left( \phi_{H_{i-1}} \right)^{-1} \right) \\ 
&=& \overline{d} \left( \phi_{H_i}, \phi_{H_{i-1}} \right) \\ 
&\le& \overline{d} \left( \phi_{H_i}, \lambda \right) + \overline{d} \left( \lambda, \phi_{H_{i-1}} \right) \\ 
&\le& \e_i + \e_{i-1} \le 2 \e_{i-1}. 
\end{eqnarray*}
For the first term, note that 
\begin{eqnarray*}
&& d_{C^0} \left( \phi_{i-1} \circ \left( \phi_{H_{i-1}} \right)^{-1} \circ \phi_{H_i}, \phi_{i-1} \right) \\ 
&& \hspace{1cm} \le d_{C^0} \left( \phi_{i-1} \circ \left( \phi_{H_{i-1}} \right)^{-1} \circ \phi_{H_i}, h \right) 
+ d_{C^0} ( h, \phi_{i-1} ) \\ 
&& \hspace{1cm} \le d_{C^0} \left( \phi_{i-1} \circ \left( \phi_{H_{i-1}} \right)^{-1}, 
h \circ \left( \phi_{H_i} \right)^{-1} \right) 
+ d_{C^0} ( \lambda, \phi_{H_{i-1}} ) \\ 
&& \hspace{1cm} \le d_{C^0} \left( \phi_{i-1} \circ \left( \phi_{H_{i-1}} \right)^{-1}, 
h \circ \left( \phi_{H_{i-1}} \right)^{-1} \right) 
+ d_{C^0} \left( h \circ \left( \phi_{H_{i-1}} \right)^{-1}, h \circ \lambda^{-1} \right) \\ 
&& \hspace{3cm} + \, d_{C^0} \left( h \circ \lambda^{-1}, h \circ \left( \phi_{H_i} \right)^{-1} \right) + \e_{i-1} \\ 
&& \hspace{1cm} \le d_{C^0} ( \phi_{i-1}, h ) + \e_{i-1} + \e_i + \e_{i-1} \\ 
&& \hspace{1cm} \le d_{C^0} ( \phi_{H_{i-1}}, \lambda ) + 3 \e_{i-1} \\ 
&& \hspace{1cm} \le 4 \e_{i-1}. 
\end{eqnarray*}
Therefore we have 
$$ \overline{d} ( \phi_{K_i}, \phi_{i-1} ) \le 4 \e_{i-1} $$ 
for all $i > 1$ as claimed. 

Now apply Lemma \ref{technicallemma} to each $K_i$ to obtain a sequence of smooth normalized Hamiltonians $L_i$ such that 
$\phi_{L_i}^0 = \phi_{K_i}^0 = \phi_{i-1}$, $\phi_{L_i}^1 = \phi_{K_i}^1 = \phi_i$, for all $i$, and 
$$ \| L_i \|_\infty < \| K_i \|_{(1,\infty)} + \e_{i-1} \le 2 \e_{i-1}, \hspace{5mm} 
\overline{d} (\phi_{L_i}, \phi_{i-1}) < \overline{d} (\phi_{K_i}, \phi_{i-1}) + \e_{i-1} \le 5 \e_{i-1} $$ 
for all $i > 1$. 

Then using Lemma \ref{inftyapproximation} to reparameterize each $L_i$ 
we obtain a smooth normalized boundary flat Hamiltonian $M_i$ such that 
$\phi_{M_i}^0 = \phi_{L_i}^0 = \phi_{i-1}$, $\phi_{M_i}^1 = \phi_{L_i}^1 = \phi_i$, for all $i$, and 
$$ \| M_i \|_\infty \le 3 \| L_i \|_\infty + \e_{i-1} \le 7 \e_{i-1}, \hspace{5mm} 
\overline{d} (\phi_{M_i}, \phi_{L_i}) < \e_{i-1} $$ 
for all $i > 1$. 
The last two inequalities can be achieved by choosing the $\e$ in Lemma \ref{inftyapproximation} to be 
$\e = \min (\e_{i-1},\frac{\e_{i-1}}{C}) > 0$, where $C$ is the constant in Lemma \ref{inftyapproximation} 
for the Hamiltonian $L_i$ (which may be different for each $i$). 
In particular, 
$$ \overline{d} (\phi_{M_i}, \phi_{i-1} ) \le \overline{d} (\phi_{M_i}, \phi_{L_i}) 
+ \overline{d} (\phi_{L_i}, \phi_{i-1} ) < 6 \e_{i-1} $$ 
for all $i > 1$. 

Finally, let $t_i = 1 - \frac{1}{2^i}$ for all $i \ge 0$. 
In particular, $0 = t_0 < t_1 < t_2 < \ldots < 1$. 
Then for $i \ge 1$ define the sequence $N_i$ of smooth normalized boundary flat Hamiltonians defined on $[t_{i-1}, t_i]$ 
by $N_i = M_i^{\zeta_{t_{i-1},t_i}}$. 
As remarked above, we have 
$$ \| N_i \|_\infty = \frac{1}{t_i - t_{i-1}} \| M_i \|_\infty = 2^i \| M_i \|_\infty < 7 \cdot 2^i \e_{i-1} $$ 
for all $i > 1$. 
By choosing $\e_i$ sufficiently small, for example $\e_{i-1} = \frac{1}{7} \frac{1}{2^i} \frac{1}{2^i}$ for $i > 1$, 
we get $\| N_i \|_\infty < \frac{1}{2^i}$, 
and since $M_i$ is just a reparameterization of $N_i$, 
$$ \overline{d} (\phi_{N_i}, \phi_{i-1} ) = \overline{d} (\phi_{M_i}, \phi_{i-1} ) < 6 \e_{i-1} < \frac{1}{2^i}. $$ 

The sequence $F_i$ of smooth normalized Hamiltonians is then defined as follows. 
Let $F_1 = N_1$ on $[0,t_1]$ and $F_1 = 0$ on $[t_1,1]$, and for $i > 1$, define 
\begin{eqnarray*}
F_i &=& F_{i-1} \hspace{1.22cm} \mbox{on } [0,t_{i-1}], \\ 
F_i &=& N_i \hspace{1.55cm} \mbox{on } [t_{i-1},t_i], \mbox{ and} \\ 
F_i &=& 0 \hspace{1.8cm} \mbox{on } [t_i,1]. 
\end{eqnarray*}
The Hamiltonians $F_i$ are indeed smooth due to boundary flatness of the $N_i$. 
We see that $ \| F_i - F_{i-1} \|_\infty = \| N_i \|_\infty < \frac{1}{2^i}$. 
In particular, $\| F_i - F_j \|_\infty \to 0$ as $i,j \to \infty$. 

It follows from the definition that $F_1$ generates a reparameterization of the path $\phi_{H_1}$, and for $i > 1$ 
the path generated by $F_i$ is equal to the one generated by $F_{i-1}$ on the interval $[0,t_{i-1}]$, 
equal to the path $\phi_{N_i}$ on the interval $[t_{i-1}, t_i]$, 
and is constant on the remaining interval $[t_i, 1]$. 
In particular, the paths $\phi_{F_i}$ are continuous, and due to the boundary flatness of the $N_i$ in fact smooth. 
Moreover the paths $\phi_{F_{i-1}}$ and $\phi_{F_i}$ agree everywhere except on the interval $[t_{i-1}, 1]$. 
Since both paths are constant on the interval $[t_i, 1]$, 
their maximum distance is achieved on the interval $[t_{i-1}, t_i]$. 
On that interval, $\phi_{F_{i-1}}$ is just the constant path $\phi_{i-1}$, 
while $\phi_{F_i}$ is the path $\phi_{N_i}$ from $\phi_{i-1}$ to $\phi_i$. 
By the above this implies that 
$$ \overline{d} (\phi_{F_{i-1}}, \phi_{F_i}) = \overline{d} (\phi_{N_i}, \phi_{i-1} ) < \frac{1}{2^i}. $$ 
In particular, $\overline{d} (\phi_{F_i}, \phi_{F_j}) \to 0$ as $i, j \to \infty$. 

That is, the sequence $(\phi_{F_i}, F_i)$ is a Cauchy sequence in $\QQ_\infty$, 
i.e., in the $L^\infty$-Hamiltonian topology. 
Since $\phi_{F_i}^1 = \phi_i \to h$ as $i \to \infty$, we conclude that $h \in Hameo_\infty \, (M,\omega)$. 
Hence $Hameo_{(1,\infty)} \, (M,\omega) \subset Hameo_\infty \, (M,\omega)$. 
The other inclusion was proved above, 
and hence we have completed the proof of Theorem \ref{maintheorem} or Theorem \ref{firsttheorem}. 
\end{proof}


\end{document}